\documentclass[11pt]{article}

\usepackage{amsmath}  % Nodig voor het gebruik van \numberwithin
\usepackage{amsthm} % nodig voor het gebruik van \begin{proof}
\usepackage{amssymb} % nodig voor \Bbb en \mathbf (wat beide hetzelfde is, maar \Bbb wordt afgeraden)

\usepackage{enumitem}
\usepackage[backref=page]{hyperref}  %

\usepackage{color} % Voor het gebruik van kleuren

\newtheorem{thm}{Theorem}
\newtheorem{inspr}[thm]{}

\newenvironment{definitie}{\begin{itemize}\item[ ]\hspace{-26pt}\bf Definition \rm }{\end{itemize}}
\newenvironment{notatie}{\begin{itemize}\item[ ]\hspace{-26pt}\bf Notation \rm }{\end{itemize}}
\newenvironment{voorbeeld}{\begin{itemize}\item[ ]\hspace{-26pt}\bf Example \rm }{\end{itemize}}
\newenvironment{stelling}{\begin{itemize}\item[ ]\hspace{-26pt}\bf Theorem \rm }{\end{itemize}}
\newenvironment{propositie}{\begin{itemize}\item[ ]\hspace{-26pt}\bf Proposition \rm }{\end{itemize}}
\newenvironment{lemma}{\begin{itemize}\item[ ]\hspace{-26pt}\bf Lemma \rm }{\end{itemize}}
\newenvironment{opmerking}{\begin{itemize}\item[ ]\hspace{-26pt}\bf Remark \rm }{\end{itemize}}
\newenvironment{voorwaarde}{\begin{itemize}\item[ ]\hspace{-26pt}\bf Condition \rm }{\end{itemize}}
\newenvironment{probleem}{\begin{itemize}\item[ ]\hspace{-26pt}\bf Problem \rm }{\end{itemize}}
\newenvironment{gevolg}{\begin{itemize}\item[ ]\hspace{-26pt}\bf Corollary \rm }{\end{itemize}}
\newenvironment{niets}{\begin{itemize}\item[ ]\hspace{-26pt}\bf   \rm }{\end{itemize}}

\newcommand{\defin}{\begin{inspr}\begin{definitie}}  %\def already defined
\newcommand{\edefin}{\end{definitie}\end{inspr}}

\newcommand{\notat}{\begin{inspr}\begin{notatie}}  %\not already defined
\newcommand{\enotat}{\end{notatie}\end{inspr}}

\newcommand{\voorb}{\begin{inspr}\begin{voorbeeld}}  %\not already defined
\newcommand{\evoorb}{\end{voorbeeld}\end{inspr}}

\newcommand{\stel}{\begin{inspr}\begin{stelling}}
\newcommand{\estel}{\end{stelling}\end{inspr}}

\newcommand{\prop}{\begin{inspr}\begin{propositie}}
\newcommand{\eprop}{\end{propositie}\end{inspr}}

\newcommand{\lem}{\begin{inspr}\begin{lemma}}
\newcommand{\elem}{\end{lemma}\end{inspr}}

\newcommand{\opm}{\begin{inspr}\begin{opmerking}}
\newcommand{\eopm}{\end{opmerking}\end{inspr}}

\newcommand{\voorw}{\begin{inspr}\begin{voorwaarde}}
\newcommand{\evoorw}{\end{voorwaarde}\end{inspr}}

\newcommand{\probl}{\begin{inspr}\begin{probleem}}
\newcommand{\eprobl}{\end{probleem}\end{inspr}}

\newcommand{\gev}{\begin{inspr}\begin{gevolg}}
\newcommand{\egev}{\end{gevolg}\end{inspr}}

\newcommand{\nul}{\begin{inspr}\begin{niets}}
\newcommand{\enul}{\end{niets}\end{inspr}}

\newcommand{\bew}{\vspace{-0.3cm}\begin{itemize}\item[ ] \bf Proof\rm: }
\newcommand{\ebew}{\hfill $\qed$ \end{itemize}}

\newcommand{\ssnl}{\vskip 3pt} % Het is noodzakelijk dat er voor de instructie vspace 
\newcommand{\nl}{\vskip 12pt} % Kunnen we dit hierin opnemen?

\newcommand{\ot}{\otimes}

\newcommand{\inv}{^{-1}}

\newcommand{\tussenen}{\qquad\quad\text{and}\qquad\quad}

\numberwithin{thm}{section}   % Zorgt ervoor dat de nummering bestaat uit ?.?
\numberwithin{equation}{section} % Zorgt ervoor dat de nummering bestaat uit ?.?

\newcommand{\keepcomment}[1]{}
\newcommand{\oldcomment}[1]{}

\begin{document}

%%%%%%%%%%%%%%%%%%%%%%%%%%%%%%%%%%%%%%%%%%%%%%%%%%%%%%%%%%%%%%%%%%%%%
%
% Abstract
%
%%%%%%%%%%%%%%%%%%%%%%%%%%%%%%%%%%%%%%%%%%%%%%%%%%%%%%%%%%%%%%%%%%%%%

\centerline{\bf \Large Single sided multiplier Hopf algebras}
\vspace{13 pt}
\centerline{\it A.\ Van Daele \rm ($^*$)}
\vspace{20 pt}
{\bf Abstract} 
\nl
Let $A$ be a non-degenerate algebra over the complex numbers and $\Delta$ a homomorphism from $A$ to the multiplier algebra $M(A\ot A)$. Consider the linear maps $T_1$ and $T_2$ from $A\ot A$ to $M(A\ot A)$ defined by
\begin{equation*}
T_1(a\ot b)=\Delta(a)(1\ot b) 
\tussenen
T_2(c\ot a)=(c\ot 1)\Delta(a).
\end{equation*}
The pair $(A,\Delta)$ is a multiplier Hopf algebra if these two maps have range in $A\ot A$ and are  bijections from $A\ot A$ to itself.
\ssnl
In \cite{VD-lsw} \emph{single sided multiplier Hopf algebras} emerge in a natural way. For this case, instead of requiring the above for the maps $T_1$ and $T_2$, we now have this property for the maps $T_1$ and $T_4$ or for $T_2$ and $T_3$ where
\begin{equation*}
T_3(a\ot b)=(1\ot b)\Delta(a)
\tussenen
T_4(c\ot a)=\Delta(a)(c\ot 1).
\end{equation*}
As it turns out, also for these single sided multiplier Hopf algebras, the existence of a unique counit and antipode can be proven. 
\ssnl
In fact, rather surprisingly, using the properties of the antipode, one can actually show that for a single sided multiplier Hopf algebra all four canonical maps are bijections from $A\ot A$ to itself. In other words, $(A,\Delta)$ is automatically a regular multiplier Hopf algebra.
\ssnl
We take the advantage of this approach to reconsider some of the known results for a regular multiplier Hopf algebra.
\nl\nl
Date: {\it 11 March 2024} 

\vskip 6cm

\hrule
\vskip 7 pt
\begin{itemize}
\item[($^*$)] Department of Mathematics, KU Leuven, Celestijnenlaan 200B,\newline
B-3001 Heverlee (Belgium). E-mail: {\tt alfons.vandaele\@kuleuven.be}
\end{itemize}

%%%%%%%%%%%%%%%%%%%%%%%%%%%%%%%%%%
%
% Introduction
%
%%%%%%%%%%%%%%%%%%%%%%%%%%%%%

\newpage
\section{\hspace{-17pt}. Introduction}
\nl
Let $A$ be an associative algebra over the field $\mathbb C$ of complex numbers. We do not require that it has an identity, but we need that the product is non-degenerate as a bilinear map. For such an algebra, we have the multiplier algebra, characterized as the largest unital algebra containing $A$ as an essential ideal. We also consider $A\ot A$ and its multiplier algebra $M(A\ot A)$. There are natural embeddings
\begin{equation*}
A\ot A\subseteq M(A)\ot M(A)\subseteq M(A\ot A)
\end{equation*}
and in general, these are strict embeddings.
\ssnl
Given a linear map $\Delta:A\to M(A\ot A)$ we consider four linear maps form $A\ot A$ to $M(A\ot A)$ defined by (see Section \ref{s:ssmhas})
\begin{align*}
T_1(a\ot b)&=\Delta(a)(1\ot b)
\tussenen
T_2 (c\ot a)=(c\ot 1)\Delta(a) \\
T_3(a\ot b)&=(1\ot b)\Delta(a)
\tussenen
T_4(c\ot a)=\Delta(a)(c\ot 1).
\end{align*}
These maps are the canonical maps associated to $\Delta$. They are called regular if they have range in $A\ot A$.
\ssnl
For a \emph{left multiplier Hopf algebra} we have a pair  $(A,\Delta)$ where $\Delta$ is a homomorphism and  the maps $T_1$ and $T_4$ are regular. Then we can require that $\Delta$ is coassociative (see Definition \ref{defin:1.11a} in Section \ref{s:ssmhas}). Finally these two maps are assumed to be bijections of $A\ot A$.  Similarly, for a \emph{right multiplier Hopf algebra}, we consider the same properties but for the maps $T_2$ and $T_3$. See Definition \ref{defin:1.3} in Section \ref{s:ssmhas}.
\ssnl
In \cite{VD-mha} multiplier Hopf algebras are defined as above, with the maps $T_1$ and $T_2$. Remark that a multiplier Hopf algebra is called regular if all four maps are regular and bijections. This means that a regular multiplier Hopf algebra is both a left and a right multiplier Hopf algebra with this terminology.
\ssnl
Such single sided multiplier Hopf algebras emerge naturally when reflecting further on the Larson Sweedler theorem as we have done in \cite{VD-lsw}. 
\ssnl
In this paper we show however that single sided multiplier Hopf algebras are in fact regular multiplier Hopf algebras. In other words, if the maps $T_1$ and $T_4$ are bijections from $A\ot A$ to itself, then this is also the case for the other maps $T_2$ and $T_3$. Similarly, when it is given that $T_2$ and $T_3$ are bijections of $A\ot A$. This is an unexpected result.
\ssnl
To obtain this result we first construct a counit and an antipode for a single sided multiplier Hopf algebra, using very much the same techniques as those used in the theory of genuine multiplier Hopf algebras. 
Nevertheless, it does not seem to be possible to derive the results for single sided multiplier Hopf algebras from those for multiplier Hopf algebras. These objects have the expected properties. 
\ssnl
In fact, along the way, one can already guess that in the end, we will arrive at a regular multiplier Hopf algebra. But it seems we need to go all the way, before we can draw that conclusion.
\ssnl
In Section \ref{s:misc} we include some proofs for known results, true for regular multiplier Hopf algebras. Doing so, we see what kind of properties of the canonical maps are really needed. Moreover, and that is in fact true for the whole paper, we use various methods that can be used in different circumstances.
\ssnl
There are natural problems that we do not treat in this paper. The main one still is to find non-trivial examples of multiplier Hopf algebras that are not regular. What we discover in this paper indicates that this may be a harder problem than first expected. 
This is still open (as far as I know). In \cite{VD-refl} we have made a step towards solving this problem by constructing some strange and non-regular coproducts on non-degenerate algebras. But this is still far from a solution to this problem.
\ssnl
We refer to the last section where we draw some conclusions and discuss possible future research. 
%\newpage%
\nl
\bf Conventions and notations \rm
\nl
We only work with algebras over $\mathbb C$ (although we believe that this is not essential and that our results are also true for algebras over other fields). The algebras are not assumed to be unital, but it is required that the product is non-degenerate. We call it a \emph{non-degenerate algebra}. Then we can consider the multiplier algebra $M(A)$ as mentioned already.
\ssnl
For an algebra $A$ we use $A^\text{op}$ for the algebra with the same underlying space but with the opposite product.
\ssnl
Occasionally we will also use the left multiplier algebra $L(A)$ and the right multiplier algebra $R(A)$. A left multiplier is a linear map $a\mapsto xa$ from $A$ to itself satisfying $x(ab)=(xa)b$ for all $a,b\in A$. A right multiplier is a linear map $a\mapsto ay$ from $A$ to itself satisfying $(ab)y=a(by)$ for all $a,b\in A$. We refer to \cite{VD-Ve1} for more details.
\ssnl
An algebra $A$ is called idempotent if any element of $A$ is the sum of products of elements of $A$. The condition is written as $A=A^2$. If the algebra has left or right local units, this is automatically fulfilled. 
\ssnl
A coproduct in this paper is a homomorphism $\Delta:A\to M(A\ot A)$ satisfying some form of coassociativity. We use $\Delta^\text{cop}$ for the new coproduct obtained from the original one by composing it with the flip map.
\ssnl
In this paper, we will sometimes use the Sweedler notation for a coproduct. This has to be done with some care because the coproduct is not assumed to have range in the tensor product $A\ot A$. One needs the regularity of the appropriate canonical maps and then to cover certain factors. It has by now become a standard technique and we will explain it further, see e.g.\ Notation \ref{notat:1.4} in Section \ref{s:ssmhas}.
\nl

\bf Basic references \rm
\nl 
For the theory of Hopf algebras we have the standard references \cite{A} and \cite{S}, see also the more recent work \cite{R}.
For the theory of multiplier Hopf algebras, the main (original) reference is \cite{VD-mha} and for the theory of multiplier Hopf algebras with integrals, sometimes called algebraic quantum groups, the main reference is \cite{VD-afw}. Weak multiplier Hopf algebras are studied in a number of papers, see \cite{VD-W1} and \cite{VD-W2} (and also \cite{VD-W0}), but we will say very little about them in this paper.
\nl
\bf Acknowledgments  \rm
\nl 
I am grateful for having the opportunity to continue my research at the KU Leuven after my retirement. Further I  like to thank my colleagues and friends at the Department of Mathematics of the University of Oslo and Trondheim (Norway), where part of this work has been written, for their hospitality and the nice working atmosphere. \oldcomment{Maybe also Trondheim and Nanjing colleagues?}

%%%%%%%%%%%%%%%%%
%
% Section 1
%
%%%%%%%%%%%%%%%%%%

\section{\hspace{-17pt}. Single sided multiplier Hopf algebras}\label{s:ssmhas}

We start with a non-degenerate algebra $A$ and a homomorphism $\Delta:A\to M(A\ot A)$. Recall the definition of the canonical maps, already formulated in the introduction.

\notat\label{notat:Tmaps}
Given a linear map $\Delta:A\to M(A\ot A)$, we consider the  maps from $A\ot A$ to $M(A\ot A)$, defined by
\begin{align*}
T_1(a\ot b)&=\Delta(a)(1\ot b)
\tussenen
T_2 (c\ot a)=(c\ot 1)\Delta(a) \\
T_3(a\ot b)&=(1\ot b)\Delta(a)
\tussenen
T_4(c\ot a)=\Delta(a)(c\ot 1)
\end{align*}
where $a,b,c\in A$.
\enotat

We call these maps the {\it canonical maps} associated to $\Delta$. If a canonical map has range in $A\ot A$, we say that it is {\it regular}. If all four of the canonical maps have range in $A\ot A$, we call $\Delta$ regular. 
\ssnl
We formulate \emph{coassociativity} of such a linear map. Here we are interested in two forms.

\defin\label{defin:1.11a} 
i) Assume that $T_1$ and $T_4$ are regular. Then we call $\Delta$ coassociative if
\begin{equation}
((\Delta\ot\iota)(\Delta(a)(1\ot b)))(c\ot 1\ot 1)
=((\iota\ot\Delta)(\Delta(a)(c\ot 1)))(1\ot 1\ot b)\label{eqn:1.3a}
\end{equation}
for all $a,b,c\in A$.
\ssnl
ii) Assume that $T_2$ and $T_3$ are regular. Then we call $\Delta$ coassociative if
\begin{equation}
(c\ot 1\ot 1)((\Delta\ot\iota)((1\ot b)\Delta(a)))
=(1\ot 1\ot b)((\iota\ot\Delta)((c\ot 1)\Delta(a)))\label{eqn1.4a}
\end{equation}
for all $a,b,c\in A$.
\edefin

There are other forms of coassociativity. The one commonly used in the theory of multiplier Hopf algebras (see \cite{VD-mha}) is valid when $T_1$ and $T_2$ are regular. The above ones are used in the study of multiplier Hopf algebroids (see \cite{T-VD}). All these different forms, when applicable, are in fact equivalent. This follows from the regularity of the product. We refer the reader to \cite{VD-refl} where the notion of coassociativity for coproducts on non-unital algebras is discussed in detail.
\ssnl
We are going to use the \emph{Sweedler notation}. It is useful for making formulas and equations more transparent. But the reader should be aware of the fact that this is just what it says, a notation for existing and well-defined expressions. The use of the Sweedler notation for a coproduct on non-unital algebras is discussed first in \cite{D-VD} and further in \cite{VD-tools}. A more recent paper about it is being prepared \cite{VD-sw}.

\notat\label{notat:1.4}
Assume that $T_1$ and $T_4$ are regular. We write
\begin{equation*}
\Delta(a)(1\ot b)=\sum_{(a)} a_{(1)}\ot a_{(2)}b
\tussenen
\Delta(a)(c\ot 1)=\sum_{(a)} a_{(1)}c\ot a_{(2)}.
\end{equation*}
We have, for all $a,b,c$, 
\begin{equation*}
(\Delta(a)(1\ot b))(c\ot 1)=(\Delta(a)(c\ot 1))(1\ot b)
\end{equation*}
and we can write this common expression as $\sum_{(a)} a_{(1)}c\ot a_{(2)}b$. By doing this, we implicitly use the equality above.
Coassociativity means 
\begin{equation*}
\sum_{(a)} \Delta(a_{(1)})(c\ot 1)\ot a_{(2)}b
=\sum_{(a)} a_{(1)}c\ot \Delta(a_{(2)})(1\ot b)
\end{equation*}
and we write this as
\begin{equation*}
\sum_{(a)} (a_{(1)}c\ot a_{(2)})\ot a_{(3)}b
=\sum_{(a)} a_{(1)}c\ot (a_{(2)}\ot a_{(3)}b).
\end{equation*}
So we can use $\sum_{(a)} a_{(1)}c\ot a_{(2)}\ot a_{(3)}b$ for this common esxpression.  Again doing this implicitly requires coassociativity.
\enotat

When using the Sweedler notation when $T_1$ and $T_4$ are regular, we always must have the right coverings by elements of $A$ on the right as we see in the above examples. We have a similar situation when $T_2$ and $T_3$ are regular. Then we must have the covering with elements on the left.
\ssnl
We illustrate this with a discussion about the \emph{pentagon equation}. It is intimately related with the above notions of coassociativity. We will also use it further (see Proposition \ref{prop:4.4}). 
\ssnl
In the next proposition we use $T$ for the canonical map $T_1$ and the leg numbering notation for $T$. This means that $T_{12}$ is the linear map from $A\ot A\ot A$ to itself with $T$ acting on the first two factors. Similarly $T_{23}$ is the one with $T$ acting on the last two factors and finally $T_{13}$ is with $T$ acting on the first and the second factor.

\prop\label{prop:1.4a}
Assume that $T_1$ and $T_4$ are regular and that $\Delta$ is coassociative as in item i ) of Definition \ref{defin:1.11a}. Then we have
$T_{23}T_{12}=T_{12}T_{13}T_{23}.$
\eprop

\bew
Let $a,b,c$ be elements of $A$. Then we have
\begin{align*}
(T_{23}T_{12})(a\ot b\ot c )
&=\sum_{(a)} T_{23}(a_{(1)}\ot a_{(2)}b\ot c)\\
&=\sum_{(a)} (a_{(1)}\ot \Delta(a_{(2)}b)(1\ot c).
\end{align*}
We multiply with an element $d$ in the first factor from the right. Then we can write
\begin{align*}
((T_{23}T_{12})(a\ot b\ot c ))(d\ot 1\ot 1)
&=\sum_{(a)} a_{(1)}d\ot \Delta(a_{(2)}b)(1\ot c)\\
&=\sum_{(a)} a_{(1)}d\ot \Delta(a_{(2)})\Delta(b)(1\ot c)\\
&=\sum_{(a),(b)} a_{(1)}d\ot a_{(2)}b_{(1)}\ot a_{(3)}b_{(2)} c.
\end{align*}
For the last step, we use coassociativity. 
\ssnl
On the other hand we have
\begin{align*}
(T_{12}T_{13}T_{23})(a\ot b\ot c )
&=\sum_{(b)}(T_{12}T_{13})(a\ot b_{(1)}\ot b_{(2)} c)\\
&=\sum_{(a),(b)} T_{12}(a_{(1)}\ot  b_{(1)}\ot a_{(2)}b_{(2)} c)\\
&=\sum_{(a),(b)} a_{(1)}\ot  a_{(2)}b_{(1)}\ot a_{(3)}b_{(2)} c.
\end{align*}
We see that 
\begin{equation*}
(T_{23}T_{12})(a\ot b\ot c )=(T_{12}T_{13}T_{23})(a\ot b\ot c )
\end{equation*}
\ebew

The reader can verify that at all places in the proof, we have the necessary coverings and the use of the Sweedler notation is justified. In fact, it is instructive to do this and get familiar with this practive.
\ssnl
It is also somewhat remarkable that one has to use in the proof that not only $T_1$ is regular, but also that $T_4$ is regular while this is not necessary for the equality.
As a matter of fact, the pentagon equation is  a form of coassociativity for the case where $\Delta$ is a homomorphism and when only $T_1$ is regular. See Definition 1.13 and Proposition 1.14 in \cite{VD-refl}.
\nl
We are now ready to define single sided multiplier Hopf algebras.

\defin\label{defin:1.3}
Let $A$ be a non-degenerate algebra and $\Delta:A\to M(A\ot A)$ a homomorphism. 
\ssnl
i) Assume that the maps $T_1$ and $T_4$ are regular and that $\Delta$ is coassociative as in item i) of Definition \ref{defin:1.11a}. Then we call $(A,\Delta)$ a \emph{left multiplier Hopf algebra} if the maps $T_1$ and $T_4$ are bijective maps of $A\ot A$ to itself.
 \ssnl
ii) Assume that the maps $T_2$ and $T_3$ are regular and that $\Delta$ is coassociative as in item ii) of Definition \ref{defin:1.11a}. Then we call $(A,\Delta)$ a \emph{right multiplier Hopf algebra} if the maps $T_2$ and $T_3$ are bijective maps of $A\ot A$ to itself.
\edefin

We clearly have that $(A,\Delta)$ is a left multiplier Hopf algebra if and only if $(A,\Delta^\text{cop})$ is a left multiplier Hopf algebra. Here $\Delta^\text{cop}$ is obtained from $\Delta$ by flipping the factors in the tensor product. On the other hand we have that $(A,\Delta)$ is a left multiplier Hopf algebra if and only if $(A^\text{op},\Delta)$ is a right multiplier Hopf algebra
\ssnl
These statements imply that results for one case give rise to similar results for the other cases. This means that we only have to prove a statement for one case, while for the other case it will already follow (or can be proven using similar arguments).

%%%%%%%%%%%%%%%%%%%%
%
% Section 2
%
%%%%%%%%%%%%%%%%%%

\section{\hspace{-17pt}. Construction and properties of the counit}\label{s:counit}

Let $(A,\Delta)$ be a left multiplier Hopf algebra as defined in \ref{defin:1.3}.

\defin
We define maps $E,E'$ from $A$ to the left multiplier algebra $L(A)$ by 
\begin{equation*}
E(a)b=mT_1\inv(a\ot b)
\tussenen
E'(a)c=m^\text{op} T_4\inv(c\ot a).
\end{equation*}
\edefin

We use $m$ for the multiplication map from $A\ot A$ to $A$ and $m^\text{op}$ for the opposite multiplication.
\ssnl
In other words, 
\begin{align*}
E(a)b&=\sum_i p_iq_i \quad\text{if}\quad a\ot b=\sum_i \Delta(p_i)(1\ot q_i)\\ 
E'(a)c&=\sum_i q_ip_i\quad\text{if}\quad  c\ot a=\sum_i \Delta(q_i)(p_i\ot 1).
\end{align*}

From these formulas, it is easy to see that indeed, $E$ and $E'$ are well-defined as left multipliers.
\ssnl
In the following proposition, we show that these left multipliers are actually scalar multiples of the identity.

\prop\label{prop:2.2}
For all $a$ we have $E(a)\in \mathbb C 1$.
\eprop
\bew
Take $a,b,c\in A$ and write $a\ot b$ as $\sum_i \Delta(p_i)(1\ot q_i)$. Then, using coassociativity,
\begin{align*}
(\Delta(a)(c\ot 1))\ot b
&=\sum_i (\Delta\ot\iota)(\Delta(p_i)(1\ot q_i))(c\ot 1\ot 1)\\
&=\sum_i (\iota\ot\Delta)(\Delta(p_i)(c\ot 1))(1\ot 1\ot q_i).
\end{align*}
From the definition of $E$ we find
\begin{align*}
\sum_{(a)} a_{(1)}c\ot E(a_{(2)})b
&=\sum_i (\Delta(p_i)(c\ot 1))(1\ot q_i)\\
&=\sum_i (\Delta(p_i)(1\ot q_i))(c\ot 1)\\
&=(a\ot b)(c\ot 1)
=ac\ot b.
\end{align*}
From the surjectivity of $T_4$ it follows that $E(a)$ is a scalar multiple of $1$ for all $a$.

\ebew

\defin
We define $\varepsilon:A\to \mathbb C$ by $E(a)b=\varepsilon(a)b$. 
\edefin
 
 We see from the proof that 
$\sum_{(a)} a_{(1)}c\ot \varepsilon(a_{(2)})b=ac\ot b$ and so 
$\sum_{(a)} a_{(1)}\varepsilon(a_{(2)})c=ac$.
We write $$\sum_{(a)} a_{(1)}\varepsilon(a_{(2)})=a.$$
%\ssnl
We now prove a similar result  for $E'$.

\prop\label{prop:2.4}
For all $a$ we have $E'(a)\in \mathbb C 1$. 
\eprop

\bew
Take $a,b,c\in A$ and write $c\ot a$ as $\sum_i \Delta(q_i)(p_i\ot 1)$. Then, using coassociativity,

\begin{align*}
c\ot ( \Delta(a)(1\ot b))
&=\sum_i(\iota\ot\Delta)(\Delta(q_i)(p_i\ot 1))(1\ot 1\ot b)\\
&=\sum_i(\Delta\ot\iota)(\Delta(q_i)(1\ot b)(p_i\ot 1\ot  1)).
\end{align*}
By the definition of $E'$ we find
\begin{align*}
\sum_{(a)} E'(a_{(1)})c\ot a_{(2)}b
&=\sum_i(\Delta(q_i)(1\ot b)(p_i\ot 1)\\
&=\sum_i\Delta(q_i)(p_i\ot 1)(1\ot b)\\
&=c\ot ab.
\end{align*}
From the surjectivity of $T_1$ it follows that $E'(a)$ is a scalar multiple of $1$ for all $a$.

\ebew

\defin
We define $\varepsilon':A\to \mathbb C$  by $E'(a)c=\varepsilon'(a)c$.
\edefin

Again we see from the proof that 
$\sum_{(a)} \varepsilon'(a_{(1)})c\ot a_{(2)}b=c\ot a$ and so  
$\sum_{(a)} \varepsilon'(a_{(1)})a_{(2)}b=ab.$
We write this as $$\sum_{(a)} \varepsilon'(a_{(1)})a_{(2)}=a.$$
\ssnl
We now show that these maps are the same.

\prop 
For all $a$ we have $\varepsilon(a)=\varepsilon'(a)$.
\eprop
\bew
For all $a$ we get from coassociativity that
\begin{align*}
(\iota\ot\varepsilon\ot\iota)\sum_{(a)} a_{(1)}c\ot a_{(2)}\ot a_{(3)}b
&=\sum_{(a)} a_{(1)}c\varepsilon (a_{(2)})\ot a_{(3)}b\\
&=\sum_{(a)} a_{(1)}c\ot a_{(2)}b.
\end{align*}
We get the same result when we apply $\iota\ot\varepsilon'\ot\iota$. In the end we find from the surjectivity of the maps $T_1$ and $T_4$ that $\varepsilon=\varepsilon'$.
\ebew

A similar argument can be used to prove the obvious uniqueness property of such a map $\varepsilon$. Then the following definition makes sense.

\defin
We call $\varepsilon$ the \emph{counit} of the left multiplier Hopf algebra $(A,\Delta)$.
\edefin

We can write the earlier formulas as in the next proposition.

\prop
The counit $\varepsilon$ of $(A,\Delta)$ satisfies for all $a$,
\begin{equation*}
\sum_{(a)} a_{(1)}\varepsilon(a_{(2)})=a
\tussenen
\sum_{(a)} \varepsilon(a_{(1)})a_{(2)}=a
\end{equation*}
as left multipliers.
\eprop

Because $\Delta$ is a homomorphism we expect the following.

\prop
The counit is a homomorphism.
\eprop

\bew
We take $a,b,c \in A$ and consider
\begin{equation*}
(\iota\ot\varepsilon)(\Delta(ab)(c\ot 1))
=abc=a(\iota\ot\varepsilon)(\Delta(b)(c\ot 1)).
\end{equation*}
We now use that $\Delta$ is a homomorphism and we replace $\Delta(b)(c\ot 1)$ by $c\ot b$. This is possible because the map $T_4$ is bijective. We find for all $a,b,c$, 
\begin{align*}
(\iota\ot\varepsilon)(\Delta(a)(c\ot b))
&=a(\iota\ot\varepsilon)(c\ot b)\\
&=ac\varepsilon(b)\\
&=(\iota\ot\varepsilon)(\Delta(a)(c\ot 1))\varepsilon(b).
\end{align*}
Finally, we replace $\Delta(a)(c\ot 1)$ by $c\ot a$ and we get
$c\varepsilon(ab)=c\varepsilon(a)\varepsilon(b)$.
\ebew

The procedure we use here to construct the counit and to find its properties is very similar to the one used in the case of a multiplier Hopf algebra (when $T_1$ and $T_2$ are assumed to be bijective). See the original paper \cite{VD-mha}.

%%%%%%%%%%%%%%%%%%%%%%%%%
%
% Section 3
%
%%%%%%%%%%%%%%%%%%%%%%%%%%%

\section{\hspace{-17pt}. The antipode and its properties}\label{s:antipode}

In this section, we show that an antipode exists with the expected properties. We first define $S$ and $S'$.

\defin
For all $a$, we define left multipliers $S(a)$ and $S'(a)$ by 
\begin{equation*}
S(a)b=(\varepsilon\ot\iota)T_1\inv(a\ot b)
\tussenen
S'(a)c=(\iota\ot\varepsilon)T_4\inv(c\ot a).
\end{equation*}
\edefin

We can rewrite these formulas as
\begin{align*}
S(a)b&=\sum_i\varepsilon(p_i)q_i \quad\text{ if }\quad a\ot b=\sum_i \Delta(p_i)(1\ot q_i)\\
S'(a)c&=\sum_i \varepsilon(q_i)p_i\quad\text{ if }\quad c\ot a=\sum_i\Delta(q_i)(p_i\ot 1).
\end{align*}

Also now, we see from these formulas that $S$ and $S'$ are well-defined linear maps from $A$ to the algebra $L(A)$ of left multipliers.
\ssnl
The following result is an immediate consequence of the definitions.

\prop
i) For all $a,b$ we have that $((\iota\ot S)\Delta(a))(1\ot b)\in A\ot A$ and  
$$m((\iota\ot S)(\Delta(a))(1\ot b))=\varepsilon(a)b.$$ 
\ssnl
ii) Similarly we have $((S'\ot\iota)\Delta(a))(c\ot 1)\in A\ot A$ and 
$$m^\text{op}((S'\ot \iota)(\Delta(a))(c\ot 1))=\varepsilon(a)c.$$
\eprop
\bew
i) We start as in the proof of Proposition \ref{prop:2.2}. Take $a,b,c\in A$ and write $a\ot b=\sum_i \Delta(p_i)(1\ot q_i)$. We get 
\begin{equation*}
(\Delta(a)(c\ot 1))\ot b
=\sum_i (\iota\ot\Delta)(\Delta(p_i)(c\ot 1))(1\ot 1\ot q_i).
\end{equation*}
By the definition of $S$ we find
\begin{align*}
(\iota\ot S)(\Delta(a)(c\ot 1))(1\ot b)
&=\sum_i((\iota\ot\varepsilon)(\Delta(p_i)(c\ot 1))(1\ot q_i)\\
&=\sum_i p_ic \ot q_i.
\end{align*}
It follows that $((\iota\ot S)\Delta(a))(1\ot b)\in A\ot A$ and that 
$$m((\iota\ot S)(\Delta(a))(1\ot b))=\varepsilon(a)b.$$
\ssnl
ii) As in the proof of Proposition \ref{prop:2.4} we start from $c\ot a=\sum_i\Delta(q_i)(p_i\ot 1)$ and write
\begin{equation*}
c\ot (\Delta(a)(1\ot b))
=\sum_i( (\Delta\ot\iota)(\Delta(q_i)(1\ot b))(p_i\ot 1\ot 1).
\end{equation*}
Apply $\omega$ on the third factor. Then
\begin{equation*}
c\ot (\iota\ot\omega)(\Delta(a)(1\ot b))
=\sum_i (\Delta((\iota\ot\omega)(\Delta(q_i)(1\ot b))))(p_i\ot 1).
\end{equation*}
From the definition of $S'$ we obtain
\begin{align*}
 S'((\iota\ot\omega)(\Delta(a)(1\ot b)))c
&=\sum_i \varepsilon((\iota\ot\omega)\Delta(q_i)(1\ot b))(p_i\ot 1)\\
&=\sum_i \omega(q_ib)p_i
\end{align*}
and hence
\begin{equation*}
 (S'\ot\iota)(\Delta(a)(1\ot b))(c\ot 1)
=\sum_i p_i\ot q_ib.
\end{equation*}
It follow that $((S'\ot \iota)\Delta(a))(c\ot 1)\in A\ot A$ and that 
$$m^\text{op}((S'\ot \iota)(\Delta(a))(c\ot 1))=\varepsilon(a)c.$$

\vskip -15pt
\ebew

We can rewrite these results using the Sweedler notation. We find
\begin{equation}
\sum_{(a)} a_{(1)} S(a_{(2)})b=\varepsilon(a)b
\tussenen
\sum_{(a)} a_{(2)}S'(a_{(1)})c=\varepsilon(a)c\label{eqn:3.1}
\end{equation}

One has to look at these formulas with some care. From the result above, we have that 
\begin{equation*}
\sum_{(a)} a_{(1)} \ot S(a_{(2)})b
\tussenen
\sum_{(a)} S'(a_{(1)})c\ot a_{(2)}
\end{equation*}
are well-defined elements in $A\ot A$. This is what we need for the formulas in Equation  (\ref{eqn:3.1}).
\ssnl
On the other hand, from the defining formulas, we also can obtain the following.

\prop
\begin{equation*}
\sum_{(a)} S(a_{(1)})a_{(2)}b=\varepsilon(a)b
\tussenen
\sum_{(a)} S'(a_{(2)})a_{(1)}c=\varepsilon(a)c
\end{equation*}
\eprop
\bew
Indeed, by the definition of $S$ we get
\begin{align*}
\sum_{(a)} S(a_{(1)})a_{(2)}b
&=(\varepsilon\ot\iota)T_1\inv(\sum_{(a)} a_{(1)}\ot a_{(2)}b))\\
&=(\varepsilon\ot \iota)(a\ot b)=\varepsilon(a)b.
\end{align*}
Similarly
\begin{align*}
\sum_{(a)} S'(a_{(2)})a_{(1)}c
&=(\iota\ot\varepsilon)T_4\inv(\sum_{(a)} a_{(1)}c\ot a_{(2)})\\
&=(\iota\ot\varepsilon)(c\ot a)=\varepsilon(a)c.
\end{align*}
\ebew

Contrary to the previous formulas, these ones are not problematic. This is because
\begin{equation*}
\sum_{(a)} a_{(1)}\ot a_{(2)}b
\tussenen
\sum_{(a)} a_{(1)}c\ot a_{(2)}
\end{equation*}
belong to $A\ot A$.
\ssnl
In the proof we also see that 
\begin{equation*}
((\iota\ot S)(\Delta(a))(1\ot b)=T_1\inv(a\ot b)
\tussenen
(S'\ot \iota)\Delta(a)(c\ot 1)=T_4\inv(c\ot a).
\end{equation*}

The information we have up to now suggest that, if $A$ is unital, then we should have a Hopf algebra with an invertible antipode $S$, the inverse being $S'$. In fact, we will find out later that we must have a regular multiplier Hopf algebra in general (see Theorem \ref{stel:3.8}). But some work still has to be done.
\ssnl
First we show that $S$ is a anti-homomorphism in the following sense.

\prop 
For all $a,b,c$ we have $S(ab)c=S(b)S(a)c$.
\eprop
\bew
i) For all $a,b,c$ we have on the one hand
\begin{equation*}
ab\ot c
=T_1 \sum_{(a),(b)} a_{(1)}b_{(1)} \ot S(a_{(2)}b_{(2)})c
\end{equation*}
while also

\begin{align*}
ab\ot c
&=(T_1 \sum_{(a)} (a_{(1)}\ot S(a_{(2)})c))(b\ot 1)\\
&=\sum_{(a)} (\Delta (a_{(1)})(1\ot S(a_{(2)})c))(b\ot 1)\\
&=\sum_{(a)} \Delta (a_{(1)})(b\ot S(a_{(2)})c)\\
&=\sum_{(a),(b)} \Delta(a_{(1)})\Delta(b_{(1)})(1\ot S(b_{(2)})S(a_{(2)})c)\\
&=\sum_{(a),(b)} \Delta (a_{(1)}b_{(1)})(1\ot S(b_{(2)})S(a_{(2)})c).
\end{align*}
One can verify that the necessary coverings are available.
Then, using that $T_1$ is injective, we get
\begin{equation*}
\sum_{(a),(b)} a_{(1)}b_{(1)} \ot S(a_{(2)}b_{(2)})c=\sum_{(a),(b)} a_{(1)}b_{(1)}\ot S(b_{(2)})S(a_{(2)})c.
\end{equation*}
We can apply $\varepsilon$ and obtain
$ S(ab)c=S(b)S(a)c$.
\ssnl
ii) In a similar way, we obtain $S'(ab)c=S'(b)S'(a)c$ for all $a,b,c$.
\ebew

We will show later that the antipode also flips the coproduct, see Proposition \ref{prop:4.4}.
\ssnl 
First we prove that $S$ and $S'$ are each others inverses. Here is the correct formulation.

\prop 
For all $a,c$ we have $S(S'(a)c)=S(c)a$ as left multipliers. Similarly we have $S'(S(a)b)=S(b)a$ for all $a,b$.
\eprop
\bew
i) Take $a,c\in A$ and apply $S$ on $\sum_{(a)} a_{(2)}S'( a_{(1)})c=\varepsilon(a)c$. Using that $S$ is a anti-homomorphism we find
\begin{equation*}
S(S'(a_{(1)})c)S( a_{(2)})=\varepsilon(a)S(c)
\end{equation*}
as left multipliers. This holds for all $a$. We can then use this formula on the left leg of $\Delta(a)(1\ot b)$ and multiply. Then we have
\begin{equation*}
S(S'(a_{(1)})c)S( a_{(2)})a_{(3)}b=S(c)\varepsilon(a_{(1)})a_{(2)}b.
\end{equation*}
For the left hand side we get eventually
\begin{equation*}
S(S'(a_{(1)})c)S( a_{(2)})a_{(3)}b=S(S'(a_{(1)})c)\varepsilon(a_{(2)})b=S(S'(a)c)b
\end{equation*}
while for the right hand side we find $S(c)ab$. This proves that $S(S'(a)c)=S(c)a$ as left multipliers.
\ssnl
ii) To prove the other statement, we just use the same argument with $S$ and $S'$ interchanged.
\ebew

As a corollary, we find that $S$ is actually a bijection of $A$ and that $S'$ is its inverse.

\prop \label{prop:3.6}
We have $S(a)\in A$ and $S'(a)\in A$ for all $a$. Consequently $S$ and $S'$ are bijections of $A$ and $S'=S\inv$.
\eprop
\bew
i) First we claim that $S(A)A=A$. 
To see this, take elements $a,b$ and consider 
\begin{equation*}
\sum_{(a)} a_{(1)} \ot S(a_{(2)})b.
\end{equation*}
This is $T_1\inv(a\ot b)$. We know that such elements span all of $A\ot A$. If we apply the counit on the first factor we find that elements of the form $S(a)b$ span all of $A$.
\ssnl
ii) From the previous result we get that $S'(S(a)c)=S(c)a$ for all $a,c$. The right hand side belongs to $A$ and so we see that also  $S'(S(a)c)\in A$. By i) we get that $S'(A)\subseteq A$. 
\ssnl
iii) A similar argument gives that also $S(A)\subseteq A$. Then the result follows.
\ebew

Compare this with the proof of Proposition 5.2 in \cite{VD-mha} where this result is proven for regular multiplier Hopf algebras.
\defin
We call $S$ the \emph{antipode} of the left multiplier Hopf algebra. 
\edefin
Proposition \ref{prop:3.6} has the following, 
rather unexpected result. 

\stel\label{stel:3.8}
 Any left multiplier Hopf algebra is a regular multiplier Hopf algebra.
 \estel
 
 \bew
 We have seen that the antipode $S$ is a bijection from $A$ to itself and that it is a anti-isomorphism. Then for all $a,b,c$ we find 
 \begin{equation*}
\sum_{(a)} a_{(1)}c \ot S(a_{(2)})S(b)=\sum_{(a)} a_{(1)}c \ot S(ba_{(2)}).
\end{equation*}
We see that 
\begin{equation*}
T_3(a\ot b)=(\iota\ot S\inv)T_1\inv(a\ot S(b)).
\end{equation*}
We see that  the map $T_3$ is a bijective map of $A\ot A$. 
\ssnl
In a similar way we get that also $T_2$ is a bijective map of $A\ot A$. Consequently, $(A,\Delta)$ is a regular multiplier Hopf algebra

 \ebew
 
Because of this result, we have to view the previous work in this paper as a result about regular multiplier Hopf algebras. It may not be very interesting as such. But nevertheless, we see another approach, using other techniques than in the original paper \cite{VD-mha}. One can see this as a new contribution to understanding the theory.
\ssnl
This is also how one should look at the next section, where we treat some known results with techniques as in the previous part of this paper.
 
 %%%%%%%%%%%%%%%%%%%%%%%%%
%
% Section 4
%
%%%%%%%%%%%%%%%%%%%%%%%%%%%

\section{\hspace{-17pt}. Miscellaneous results revisited}\label{s:misc}

In this section, we will prove some properties of regular multiplier Hopf algebras that are scattered in the literature. We not only provide simpler arguments, but we also indicate what  conditions are really needed to get these results. This material has a natural place in this note.
\ssnl
We assume in this section that $A$ is a non-degenerate algebra and that $\Delta:A\to M(A\ot A)$ is just a linear map to begin with. We consider the canonical maps $T_1$, $T_2$, $T_3$ and $T_4$ as in Notation \ref{notat:Tmaps}.
\ssnl
First observe that $\Delta$ is automatically full if the appropriate canonical maps are regular and bijective. If e.g.\ $T_1$ is regular and bijective, we have that the span of elements of the form $(\iota\ot\omega)(\Delta(a)(1\ot b))$ is all of $A$. Hence the left leg of $\Delta$ is all of $A$. For the property of fullness, we refer to the recent paper \cite{VD-refl}. 
\ssnl
The following is a little less trivial, again see \cite{VD-refl}. 

 \prop
 Assume that $T_1$ and $T_4$ are regular and that the range of $T_1$ is all of $A\ot A$. Then $A$ is idempotent.
  \eprop
 
 \bew
 Assume that $\omega$ is a linear functional on $A$ that is $0$ on $A^2$. Then, for all $a,b,c$ we have 
 \begin{equation*}
(\iota\ot\omega)(\Delta(a)(c\ot b))=0.
\end{equation*}
We use that $\Delta(a)(c\ot 1)\in A\ot A$. Now  because also $\Delta(a)(1\ot b)\in A\ot A$ we can cancel $c$ and obtain that
 \begin{equation*}
(\iota\ot\omega)(\Delta(a)(1\ot b))=0.
\end{equation*}
Finally, because $T_1$ is assumed to be surjective, we must have that $\omega=0$. Consequently $A^2=A$.
 \ebew
 
It follows that the underlying algebra of a regular multiplier Hopf algebra is idempotent. It is clear that there are many similar conditions that can be used to get the same result. If e.g.\ $T_1$ and $T_2$ are regular and $T_1$ surjective, we also get that $A$ is idempotent. Now it follows that the underlying algebra of any multiplier Hopf algebra is idempotent. Remark in passing that such an argument can not be used for weak multiplier Hopf algebras. This is one of the reasons why, for a weak multiplier Hopf algebra, it is assumed from the very beginning that the underlying algebra is idempotent, see \cite{VD-W1}.
\ssnl
A simple consequence is the following.

\prop
Assume that $T_1$ and $T_4$ are regular. If the range of $T_1$ is all of $A\ot A$, then $\Delta(A)(A\ot A)=A\ot A$.
\eprop 
\bew
We have, using that the range of $T_1$ is $A\ot A$,
\begin{equation*}
\Delta(A)(A\ot A)=\Delta(A)(1\ot A)(A\ot 1)=A^2\ot A.
\end{equation*}
From the previous property we get $A^2=A$ and this proves the result.
\ebew

We cannot conclude from this that $\Delta$ is non-degenerate. We also need $(A\ot A)\Delta(A)=A\ot A$. For this  we would need regularity of the other canonical maps. Fortunately, for a regular multiplier Hopf algebra, we all have these results and hence, the coproduct is non-degenerate in that case.
\ssnl
Next we look at the existence of local units in $A$. This is much stronger than $A^2=A$.

\prop
Assume that $(A,\Delta)$ is a regular multiplier Hopf algebra. 
Then the algebra $A$  has  left local units.
\eprop
\bew
i) Take $b\in A$ and assume that $\omega$ is a linear functional that is $0$ on $Ab$. 
Take $a\in A$ and write
\begin{equation*}
a\ot b=\sum_i \Delta(p_i)(1\ot q_i).
\end{equation*}
For all $c, p\in A$ we have
\begin{equation*}
(\iota\ot\omega)(\sum_{(a)} a_{(1)}c\ot pS(a_{(2)})b)=0.
\end{equation*}
We use that $\sum_{(a)} a_{(1)}c\ot a_{(2)}\in A\ot  A$ and that $S$ maps into $A$. Further, because also $\sum_{(a)} a_{(1)}\ot S(a_{(2)})b\in A\ot  A$, we can cancel $c$ and obtain that
\begin{equation*}
(\iota\ot\omega)(\sum_{(a)} a_{(1)}\ot pS(a_{(2)})b)=0.
\end{equation*}
Finally, we use again that $\sum_{(a)} a_{(1)}\ot S(a_{(2)})b\in A\ot  A$ and we can replace $p$ by $a_{(1)}$
to get 
\begin{equation*}
\omega(\sum_{(a)} a_{(1)}S(a_{(2)})b)=0.
\end{equation*}
This implies that $\varepsilon(a)\omega(b)=0.$ As a consequence we find $b\in Ab$. Then $A$ has left local units.
\ssnl
ii) In a similar way, or by using that $S$ is an anti-isomorphism of $A$, we get that $A$ has right local units. Consequently it has two-sided local units by an argument in \cite{Ve}. 
\ebew

The proof also works for non-regular multiplier Hopf algebras. In that case we use that $S$ maps $A$ into the multiplier algebra. See e.g.\ \cite{VD-Ve1}.
\keepcomment{Include more references}{\ssnl}

We see that in order to obtain this result, we need more properties of the canonical maps. It would be worthwhile to find out what is precisely needed.

\oldcomment{We also need to make remarks about this. One can try to find out what properties are really necessary to obtain this result. }{}
\ssnl
We now show that the antipode flips the coproduct.

\prop \label{prop:4.4}
For all $a,b,c$  we have
\begin{equation*}
\Delta(S(a)b)(1\ot c)=\sum_{(a)} (S(a_{(2)}) \ot S(a_{(1)}))\Delta(b)(1\ot c).
\end{equation*}
\eprop

\bew
We denote in this proof the canonical map $T_1$ by $T$ and use the leg-numbering notation for the map $T$ as in Proposition \ref{prop:1.4a}. From the pentagon equation for $T_1$ we find
\begin{equation*}
T_{13}\inv T_{12}\inv T_{23}=T_{23}T_{12}\inv. 
\end{equation*}
We apply this on $a\ot b\ot c$. For the left hand side we obtain
\begin{align*}
(T_{13}\inv T_{12}\inv T_{23})(a\ot b\ot c)
&=\sum_{(b)} T_{13}\inv T_{12}\inv  (a\ot b_{(1)}\ot b_{(2)} c ) )\\
&=\sum_{(a),(b)} T_{13}\inv (a_{(1)} \ot S(a_{(2)})b_{(1)}\ot b_{(2)} c )\\
&=\sum_{(a),(b)}  (a_{(1)} \ot S(a_{(3)})b_{(1)}\ot S(a_{(2)})b_{(2)} c.
\end{align*}
For the right hand side we find
\begin{align*}
(T_{23}T_{12}\inv)(a\ot b\ot c)
&=\sum_{(a)}T_{23} (a_{(1)}\ot S(a_{(2)})b\ot c)\\
&=\sum_{(a)}a_{(1)}\ot \Delta(S(a_{(2)}b)(1 \ot c).
\end{align*}
These two expressions are equal and if we apply the counit on the first factor we obtain
\begin{equation*}
\sum_{(a)}  (S(a_{(2)})\ot S(a_{(1)})\Delta(b)(1\ot c)=\Delta(S(a)b)(1\ot c).
\end{equation*}
\ebew

One can verify that the necessary coverings are present so that the use of the Sweedler notation is justified.
\ssnl
It is also interesting to compare this proof with the one given in the original paper \cite{VD-mha} of this result.

%%%%%%%%%%%%%%%%%%%
%
% Section 5
%
%%%%%%%%%%%%%%%

\section{\hspace{-17pt}. Conclusions and further research}\label{s:concl}
\ssnl
In this paper, we have developed a theory of single-sided multiplier Hopf algebras. The different steps are very much like the ones for genuine multiplier Hopf algebras as found in the original paper \cite{VD-mha}. So, there is a counit and an antipode with the expected properties.
\ssnl
Still, it does not seem possible to derive these results as consequences of the results for multiplier Hopf algebras.
\ssnl
On the other hand, it turns out that any single-sided multiplier Hopf algebra is in fact a regular multiplier Hopf algebra. This is somewhat remarkable as it shows that properties of two of the canonical maps turn out to be sufficient to get properties of the other ones.
\ssnl
In this paper we make use of the Sweedler notation from the very beginning, making the results and arguments easier to follow that in the original paper. We explain why this is justified. We also take advantage of this approach to give proofs of some known results for regular multiplier Hopf algebras.
\ssnl
We are convinced that this note is a contribution to understanding and working with this theory.
\ssnl
One could try a similar approach to the theory of weak multiplier Hopf algebras as treated in \cite{VD-W0} and \cite{VD-W1}.
\ssnl
One also could review the theory of multiplier algebroids in \cite{T-VD}.
\ssnl
Finally, we again emphasize the problem of finding non-trivial examples of (weak) multiplier Hopf algebras that are not regular. In particular one has to find cases where the antipode does not map into $A$ but only into the multiplier algebra $M(A)$. This is related with finding examples of non-regular separability idempotent (see \cite{VD-si.v1, VD-si.v2}.  We have an attempt in this direction in \cite{VD-refl}.

%%%%%%%%%%%%%%%%%%%%%%%%%%%%%%%%%%%%%%%%
%
% References
%
%%%%%%%%%%%%%%%%%%%%%%%%%%%%%%%%%%%%%%%%% 

\end{document}